\documentclass{amsart}
\usepackage{amssymb,latexsym,amsmath,amsfonts,hhline}
\usepackage{pdfsync,comment}
\usepackage[usenames]{color} 

\input xy
\xyoption{all}

\newcommand{\cal}{\mathcal}

\newtheorem{formula}{}[section]
\newtheorem{definition}[formula]{Definition}
\newtheorem{corollary}[formula]{Corollary}
\newtheorem{proposition}[formula]{Proposition}
\newtheorem{remark}[formula]{Remark}
\newtheorem{lemma}[formula]{Lemma}
\newtheorem{theorem}[formula]{Theorem}

\def\thrm{\begin{theorem}}
\def\thrml#1{\begin{theorem}\label{#1}}
\def\ethrm{\end{theorem}}
\def\rmrk{\begin{remark}}
\def\rmrkl#1{\begin{remark}\label{#1}}
\def\ermrk{\end{remark}}
\def\dfntn{\begin{definition}}
\def\dfntnl#1{\begin{definition}\label{#1}}
\def\edfntn{\end{definition}}
\def\nmrt{\begin{enumerate}}
\def\enmrt{\end{enumerate}}
\def\tm#1{\item[{\rm (#1)}]}
\def\qtn{\begin{equation}}
\def\qtnl#1{\begin{equation}\label{#1}}
\def\eqtn{\end{equation}}
\def\lmm{\begin{lemma}}
\def\lmml#1{\begin{lemma}\label{#1}}
\def\elmm{\end{lemma}}
\def\crllr{\begin{corollary}}
\def\crllrl#1{\begin{corollary}\label{#1}}
\def\ecrllr{\end{corollary}}
\def\prpstn{\begin{proposition}}
\def\prpstnl#1{\begin{proposition}\label{#1}}
\def\eprpstn{\end{proposition}}
\def\css{\begin{cases}}
\def\ecss{\end{cases}}

\def\proof{\noindent{\bf Proof}.\ }

\def\cG{{\cal G}}

\def\cL{{\cal L}}

\def\cT{{\cal T}}

\def\cX{{\cal X}}
\def\cY{{\cal Y}}

\def\fG{{\mathfrak G}}
\def\fF{{\mathfrak F}}

\def\fS{{\mathfrak S}}

\def\mF{{\mathbb F}}

\DeclareMathOperator{\aut}{Aut}

\DeclareMathOperator{\diag}{Diag}
\DeclareMathOperator{\dom}{Dom}

\DeclareMathOperator{\GL}{GL}

\DeclareMathOperator{\Hom}{Hom}
\DeclareMathOperator{\id}{id}
\DeclareMathOperator{\im}{im}

\DeclareMathOperator{\iso}{Iso}

\DeclareMathOperator{\orb}{Orb}

\DeclareMathOperator{\sym}{Sym}

\def\eprf{\hfill$\square$}

\def\qaq{\quad\text{and}\quad}

\def\ov{\overline}
\newcommand{\grp}[1]{\langle {#1}\rangle}

\def\wh{\widehat}
\def\wt{\widetilde}

\def\mattwo#1#2#3#4{\begin{pmatrix} #1 & #2 \\ #3 & #4 \\ \end{pmatrix} }

\newcommand{\forme}[1]{}

\begin{document}
\title{Schurity and separability\\ of quasiregular coherent configurations}

\author{Mitsugu Hirasaka}
\address{Department of Mathematics, Pusan National University, Jang-jeon dong, Busan, Republic of Korea}
\email{hirasaka@pusan.ac.kr}
\thanks{This research was supported by Basic Science Research Program through the National Research Foundation of Korea(NRF) funded by the Ministry of  Science, ICT \& Future Planning(NRF-2016R1A2B4013474).}
\author{Kijung Kim}
\address{Department of Mathematics, Pusan National University, Jang-jeon dong, Busan, Republic of Korea}
\email{knukkj@pusan.ac.kr}
\thanks{The second author was supported by Basic Science Research Program through the National Research Foundation of Korea funded by the Ministry of Education (NRF-2017R1D1A3B03031349)}
\author{Ilia Ponomarenko}
\address{St.Petersburg Department of the Steklov Mathematical Institute, St.Petersburg, Russia}
\email{inp@pdmi.ras.ru}
\thanks{The work of the third author was supported by the RAS Program of Fundamental Research ``Modern Problems of Theoretical Mathematics''. }

\date{December 31, 2017}

\begin{abstract}
A permutation group is said to be quasiregular if every its transitive constituent is regular, and a quasiregular coherent configuration can be thought as a combinatorial analog of such a group: the transitive constituents are replaced by the homogeneous components. In this paper, we are interested in the question when the configuration is schurian, i.e., formed by the orbitals of a permutation group, or/and separable, i.e., uniquely determined by the intersection numbers. In these terms, an old result of Hanna Neumann	is, in a sense, dual to the statement that the quasiregular coherent configurations with cyclic homogeneous components are schurian. In the present paper, we (a) establish the duality in a precise form and (b) generalize the latter result by proving that a quasiregular coherent configuration is schurian and separable if the groups associated with homogeneous components have distributive lattices of normal subgroups.
\end{abstract}

\maketitle

\section{Introduction}

A permutation group is said to be {\it quasiregular} if every its transitive constituent is regular~\cite[p.~53]{DP}. This concept has an obvious direct analog in the theory of ``groups without groups'': a coherent configuration is said to be quasiregular if each its homogeneous component is regular (for exact definitions concerning coherent configurations, see Section~\ref{200817a}). In fact, there is a one-to-one correspondence between the regular groups (which are exactly transitive quasiregular groups) and regular coherent configurations~\cite{Z5}. However, no such a correspondence exists for intransitive quasiregular groups: the reason is that not every quasiregular coherent configuration is {\it schurian}, i.e., formed by the orbitals of a permutation group. The first example of a non-schurian quasiregular coherent configuration was constructed by Sergei Evdokimov in the end of 1990s (but was never published). One of motivations for the present paper is to find conditions for a quasiregular coherent configuration to be schurian.\medskip

Another motivation comes from an old paper of Hanna Neumann \cite{N50},\footnote{The third author is expressed gratitude to Mikhail Muzychuk, who drew his attention to paper~\cite{N50}.} where she studied the amalgams of finite cyclic groups. Roughly speaking, the question is whether a
family $\fG$ of groups with prescribed pairwise intersections {\it admits an amalgam}, i.e., can be isomorphically embedded to a certain group. The main result of \cite{N50} states that if $\fG$ consists of finite cyclic groups, then this is true if some natural necessary conditions concerning the pairwise intersections are satisfied. These conditions have sense not only for cyclic groups and for any family $\fG$ satisfying these conditions, we use term {\it amalgam configuration} based on~$\fG$ (the exact definition is given in Section~\ref{210817a}). Not every amalgam configuration based on a family $\fG$ of abelian groups admits an amalgam. A reason for this is revealed in the theorem below establishing a close relationship between the amalgam configurations based on~$\fG$ and quasiregular coherent configuration of type $\fG$, i.e., those the homogeneous components of which are exactly the groups from~$\fG$.

\thrml{241216g}
Let $\fG$ be a family of finite abelian groups. Then the quasiregular coherent configurations of type $\fG$ are in one-to-one correspondence with the amalgam configurations based on~$\fG$. Moreover,  a coherent configuration of type $\fG$ is schurian if and only if the corresponding amalgam configuration admits an amalgam.
\ethrm

The proof of Theorem~\ref{241216g} is presented in Section~\ref{220817c} and uses a quite general concept of a system of linked sections based on a family $\fG$ of groups (see Section~\ref{210817a}). The special cases of this concept are the amalgam configurations  and quasiregular coherent configurations; the fact that the latter are in one-to-one correspondence with the systems of linked quotients is proved in Section~\ref{040816a}. The main idea of the whole proof is to show that if the family $\fG$ consists of abelian groups, then these two special cases are dual each to other in the sense of the duality theory of abelian groups.\medskip 

From Theorem~\ref{241216g}, it immediately follows that the above mentioned result of Hanna Neumann is equivalent to the fact that any quasiregular coherent configuration 
such that all its homogeneous components are cyclic, is schurian.
Since the subgroup lattice of a cyclic group is distributive, Theorem~\ref{130617a} below generalizes the latter statement to a   broader class of groups containing all cyclic groups. Before we state this result in a precise form, we discuss one more concept making this theorem even more stronger.\medskip

Let $K\le\sym(\Omega)$, and let $\cX$ be the coherent configuration formed by the orbitals of~$K$. A natural invariant of $\cX$ and hence of~$K$ is given by the tensor of intersection numbers: they can be thought as the structure constants of the the centralizer algebra of~$K$ with respect to the linear basis formed by the adjacency matrices of the orbitals of~$K$. The question going back to D.~Higman \cite{H70} is when this tensor determines the group~$K$ up to permutation isomorphism. In the language of the coherent configurations, the positive answer exactly means that $\cX$ is {\it separable}. Now, we are ready to state the second main result of the present paper.

\thrml{130617a}
Let $\fG$ be a family of groups with  distributive lattices  of normal subgroups. Then any quasiregular coherent configuration of type $\fG$ is schurian and separable.
\ethrm

We deduce this result from Theorem~\ref{distributive} proved in Section~\ref{220817y} for even more broader class of coherent configurations. The idea of the proof is close to that used in P.~H.~Zieschang's paper~\cite{Zies}, where the schurity part of Theorem~\ref{130617a} was proved for a class of homogeneous coherent configurations. Note that these configurations are nothing else than homogeneous algebraic fusions of quasiregular coherent configuration of type $\fG$ satisfying the condition of Theorem~\ref{130617a} and such that all the groups in $\fG$ are isomorphic.\medskip

In Section~\ref{220817u}, we study the schurity and separability of all quasiregular coherent configurations with small number of homogeneous components. The main result here is presented by the following theorem.

\thrml{120817d}
A quasiregular coherent configuration with at most three homogeneous components is schurian and separable.
\ethrm

In the case of commutative homogeneous components, the schurity part of Theorem~\ref{120817d} was proved in~\cite{F89}\footnote{Theorem~5.6 in~\cite{F89} states that any  quasiregular coherent configuration with commutative homogeneous components is separable. However, this is not true: an infinite family of counterexamples was constructed in~\cite[Section~5]{EP99}.}. In Section~\ref{220817u}, we also construct an infinite family of non-schurian quasiregular coherent configurations with four commutative homogeneous components. The question of separability of these coherent configurations remains open.\medskip

The class of meta-thin association schemes was introduced and studied by the first author in~\cite{H05}. From \cite[Theorem~2.1]{EP10}, it follows that any such scheme is an algebraic fusion of a quasiregular coherent configuration with isomorphic homogeneous components. The number of them equals the index of the thin residue of the original association scheme. Thus, as a corollary of Theorem~\ref{120817d} we obtain the following statement also proved in Section~\ref{220817u}.

\crllrl{120817e}
Any meta-thin scheme with thin residue of index at most three is schurian and separable.
\ecrllr

{\bf Notation.}

Throughout the paper, $\Omega$ denotes a finite set.

The diagonal of the Cartesian product $\Omega\times\Omega$ is denoted by~$1_\Omega$.

For a relation $r\subseteq\Omega\times\Omega$, we set $r^*=\{(\beta,\alpha):\ (\alpha,\beta)\in r\}$ and $\alpha r=\{\beta\in\Omega:\ (\alpha,\beta)\in r\}$ for all $\alpha\in\Omega$.

For a relation $r\subseteq\Omega\times\Omega$ and sets $\Delta,\Gamma\subseteq\Omega$, we set $s_{\Delta,\Gamma}=s\cap(\Delta\times\Gamma)$. If $S$ is a set of relations, we put $S_{\Delta,\Gamma}=\{s_{\Delta,\Gamma}:\ s\in S\}$ and denote $S_{\Delta,\Delta}$ by $S_\Delta$.

For relations $r,s\subseteq \Omega\times\Omega$, we set
$r\cdot s=\{(\alpha,\beta):\ (\alpha,\gamma)\in r,\ (\gamma,\beta)\in s$
for some $\gamma\in\Omega\}$. If $S$ and $T$ are sets of relations, we set
$S\cdot T=\{s\cdot t:\ s\in S,\, t\in T\}$.

For a set $S$ of relations on~$\Omega$, we denote by $S^\cup$ the set of all unions of the elements of $S$, and put $S^*=\{s^*:\ s\in S\}$
and $\alpha S=\cup_{s\in S}\alpha s$, where $\alpha\in\Omega$.

\section{Coherent configurations}\label{200817a}

In our presentation of coherent configurations, we follow papers~\cite{MP,EP09} and monograph~\cite{Z5}. All the facts we use, can be found in these sources and references therein.

\subsection{Definitions.}
A pair $\cX=(\Omega,S)$, where $\Omega$ is a finite set and $S$ is a partition of $\Omega\times\Omega$, is called a {\it coherent configuration} on $\Omega$ if $1_\Omega\in S^\cup$, $S^*=S$, and if given $r,s,t\in S$, the number
$$
c_{rs}^t=|\alpha r\cap \beta s^*|
$$
does not depend on a choice of $(\alpha,\beta)\in t$. The elements of $\Omega$, $S$, $S^\cup$, and the numbers~$c_{rs}^t$ are called the {\it points}, {\it basis relations}, {\it relations} and {\it intersection numbers} of~$\cX$, respectively. The numbers $|\Omega|$ and $|S|$ are called the {\it degree} and the {\it rank} of~$\cX$. A unique basic relation containing a pair $(\alpha,\beta)\in\Omega\times\Omega$ is denoted by $r(\alpha,\beta)$. Since the mapping $r:\Omega\times \Omega \to S$ depends only on $\cX$, it should be denoted by $r_\cX$, but we usually omit the subindex if this does not lead to confusion.\medskip

The set $S^\cup$ contains the relation $r\cdot s$ for all $r,s\in S^\cup$. It follows that this relation is the union (possibly empty) of basis relations of~$\cX$; the set of these relations is called the {\it complex product} of $r$ and $s$ and denoted by $rs$. Thus 
$$
rs\subseteq S,\qquad r,s\in S.
$$
In what follows, for any $X,Y\subseteq S$, we denote by $XY$ the union of all sets $rs$ with $r\in X$ and $s\in Y$. Obviously, $(XY)Z=X(YZ)$ for all~$X,Y,Z\subseteq S$. A nonempty subset $X$ of $S$ is said to be {\it closed} if $XX^\ast \subseteq X$.

\subsection{Fibers and homogeneity.}
A set $\Delta\subseteq\Omega$ is called a {\it fiber} of the coherent configuration~$\cX$ if $1_\Delta\in S$; the set of all fibers is denoted by $F=F(\cX)$. The point set~$\Omega$ is a disjoint union of fibers. If $\Delta$ is a union of fibers, then the pair 
$$
\cX_\Delta=(\Delta,S_\Delta)
$$
is a coherent configuration, called the {\it restriction} of~$\cX$ to~$\Delta$; it is called a {\it homogeneous component} of~$\cX$ if $\Delta\in F$. The coherent configuration is said to be {\it homogeneous} or {\it association scheme} or {\it scheme} if it has exactly one homogeneous component, or equivalently if $1_\Omega\in S$.\medskip

For any basic relation $s\in S$, there exist uniquely determined fibers $\Delta,\Gamma$ such that $s\subseteq\Delta\times\Gamma$; in particular, the union
$$
S=\bigcup_{\Delta,\Gamma\in F}S_{\Delta,\Gamma}
$$
is disjoint. The number $|\delta s|$ with $\delta\in\Delta$, equals the intersection number $c_{ss^*}^{1_\Delta}$, and hence does not depend on the choice of the point~$\delta$. It is called the {\it valency} of $s$ and denoted by $n_s$. 
For each fiber $\Delta$ and $s\in S_\Delta$ we have $n_s=n_{s^\ast}$.

\subsection{Isomorphisms and schurity.}
Two coherent configurations are called {\it isomorphic} if there exists a bijection between their point sets that induces the bijection between their sets of basis relations. Each such bijection is called an {\it isomorphism} between these two configurations. The group of all isomorphisms of a  coherent configuration $\cX=(\Omega,S)$ to itself
contains a normal subgroup
$$
\aut(\cX)=\{f\in\sym(\Omega):\ s^f=s,\ s\in S\}
$$
called the {\it automorphism group} of~$\cX$, where
$s^f=\{(\alpha^f,\beta^f):\ (\alpha,\beta)\in s\}$.\medskip

Conversely, let $K\le\sym(\Omega)$ be a permutation group, and let $S=\orb(K,\Omega^2)$. Then, $\cX=(\Omega,S)$ is a coherent configuration;
we say that $\cX$ {\it is  associated} with~$K$. A coherent configuration on $\Omega$ is said to be {\it schurian} if it is associated with some permutation group on~$\Omega$. A coherent configuration~$\cX$ is schurian if and only if it is associated with the group~$\aut(\cX)$.

\subsection{Algebraic isomorphisms and separability.}
Let $\cX=(\Omega,S)$ and $\cX'=(\Omega',S')$ be coherent configurations. A bijection $\varphi:S\to S',\ r\mapsto r'$ is called an {\it algebraic isomorphism} from~$\cX$ onto~$\cX'$ if
\qtnl{f041103p1}
c_{rs}^t=c_{r's'}^{t'},\qquad r,s,t\in S.
\eqtn
In this case, $\cX$ and $\cX'$ are said to be {\it algebraically isomorphic}. Each isomorphism~$f$ from~$\cX$ onto~$\cX'$ induces an algebraic isomorphism $\varphi_f:r\mapsto r^f$ between these configurations. The set of all isomorphisms inducing the algebraic isomorphism~$\varphi$ is denoted by $\iso(\cX,\cX',\varphi)$. In particular,
\qtnl{190316b}
\iso(\cX,\cX,\id_S)=\aut(\cX),
\eqtn
where $\id_S$ is the identity mapping on $S$. A coherent configuration~$\cX$ is said to be {\it separable} if for any algebraic isomorphism~$\varphi:S\to S'$, the set $\iso(\cX,\cX',\varphi)$ is not empty.\medskip

The algebraic isomorphism $\varphi$ induces a bijection from $S^\cup$
onto $(S')^\cup$: the union $r\cup s\cup\cdots$ of basis relations
of $\cX$ is taken to $r'\cup s'\cup\cdots$. This bijection is also
denoted by $\varphi$. One can see that $\varphi$ preserves the reflexive
basis relations on all fibers. This extends $\varphi$ to a bijection
$F(\cX)\to F(\cX')$ so that $(1_{\Delta^{}})'=1_{\Delta'}$.

\subsection{Faithful maps.}
Let $\cX=(\Omega,S)$ and $\cX'=(\Omega',S')$ be coherent configurations, and let $\varphi:S\to S'$ be 
an algebraic isomorphism. A bijection $f$ from a subset of $\Omega$ to a subset of $\Omega'$ is said to be $\varphi$-\textit{faithful} if
$$
r(\alpha,\beta)^\varphi=r'(\alpha^f,\beta^f)\quad
\text{for all}\ \,\alpha,\beta\in\dom(f),
$$
where $\dom(f)$ is the domain of~$f$. Note that if this domain is a singleton $\{\alpha\}$ and $\Delta$ is the fiber of~$\cX$ containing~$\alpha$, then the fiber of~$\cX'$ containing $\alpha^f$ is equal to $\Delta^\varphi$. Clearly, if $f$ is a $\varphi$-faithful map, then the restriction of $f$ to any subset of~$\dom(f)$ is also $\varphi$-faithful.\medskip

A $\varphi$-faithful map $f$ is said to be {\it $\varphi$-extendable} to a point $\gamma\in\Omega$ if there exists a $\varphi$-faithful map with domain $\dom(f)\,\cup\,\{\gamma\}$, or, equivalently, if 
\qtnl{100617u}
\bigcap_{\alpha\in\dom(f)}\alpha^fr(\alpha,\gamma)^\varphi\ne\varnothing.
\eqtn
A $\varphi$-faithful map which is $\varphi$-extendable  to every point of~$\Omega$, is said to be $\varphi$-extendable. From the definitions of  coherent configurations and algebraic isomorphisms,
it follows that every $\varphi$-faithful map $f$ with $|\dom(f)|\le 2$ is $\varphi$-extendable. In these terms, one can give a sufficient condition for schurity and separability of a coherent configuration.

\lmml{seps}
Let $\cX=(\Omega,S)$ be a coherent configuration. Then
\nmrt
\tm{1} $\cX$ is schurian if every $\varphi$-faithful map with $\varphi=\id_S$ is $\varphi$-extendable,
\tm{2} $\cX$ is separable if for every algebraic isomorphism $\varphi$ from $\cX$ onto another coherent configuration, each $\varphi$-faithful map is $\varphi$-extendable.
\enmrt
\elmm
\proof
Let $\varphi:\cX\to\cX'$ be an algebraic isomorphism. Assume that each $\varphi$-faithful map is $\varphi$-extendable. Then for any pairs $(\alpha,\beta)\in\Omega\times\Omega$ and  $(\alpha',\beta')\in\Omega'\times\Omega'$ such that 
$$
r(\alpha,\beta)^\varphi=r'(\alpha',\beta'), 
$$
there exists an isomorphism $f\in\iso(\cX,\cX',\varphi)$ taking $\alpha$ to $\alpha'$, and $\beta$ to~$\beta'$. Now if $\cX=\cX'$ and $\varphi$ is the identity algebraic isomorphism, then $f\in\aut(\cX)$ and the coherent configuration $\cX$ is schurian. This proves statement~(1). Also, if $\varphi$ runs over all algebraic isomorphisms from $\cX$, then $\cX$ is separable, which proves statement~(2).\eprf

\crllr
Let $\cX$ be a coherent configuration. Suppose that for every algebraic isomorphism $\varphi$ onto another coherent configuration, each $\varphi$-faithful map is $\varphi$-extendable.
Then $\cX$ is schurian and separable.
\ecrllr

\subsection{Quasiregular coherent configurations.}
A relation $s\subseteq\Omega\times\Omega$ is said to be {\it thin} if 
$|\alpha s|\le 1$ and $|\beta s^*|\le 1$ for all $\alpha,\beta\in\Omega$, or equivalently, $s$ is a bijection from $\dom(s)$ onto $\dom(s^*)$. In the latter sense, $s^*$ is the bijection inverse to~$s$. One can easily see that if $\cX=(\Omega,S)$ is a coherent configuration and $s\in S$, then 
$n_{s^{}}=n_{s^*}=1$ and $s\cdot r,r\cdot s\in S$ for all $r\in S$ unless $s\cdot r$ or $s\cdot r$ are empty.

\lmml{sss}
Let $\cX=(\Omega,S)$ be a coherent configuration and $s\in S$. Assume that
the set $ss^*$ consists of thin relations. Then $ss^* s=\{s\}$.
\elmm
\proof 	Let $t\in ss^\ast$. Then $t$ is a thin relation by the assumption and hence  $ts$ is a singleton. Since also $s\in ts$, it follows  that $ts=\{s\}$. This is true for all $t\in ss^\ast$ and hence $ss^* s\subseteq \{s\}$. However, $1_\Delta\in ss^*$, where $\Delta$ is a unique fiber of $\cX$ such that $s\subseteq\Delta\times\Omega$. Thus, $s=1_\Delta\cdot s\in ss^*s$ and so $ss^*s=\{s\}$.\eprf\medskip

A coherent configuration $\cX$ is said to be {\it semiregular} if every its basis relation is thin. If, in addition, it is homogeneous, then we say that $\cX$ is a {\it regular scheme} (or thin scheme in the terminology of~\cite{Zies}). Finally, $\cX$ is said to be {\it quasiregular} if each homogeneous component of $\cX$ is regular. Note that if $K$ is a permutation group and $\cX$ is the coherent configuration associated with~$K$, then $\cX$ is semiregular (respectively, regular, quasiregular) if and only if the group~$K$ is semiregular (respectively, regular, quasiregular). It should also be noted that the basis relations of a regular scheme form a group with respect to the product $\cdot$ of relations.

\section{Systems of linked sections}\label{210817a}

In this section, we introduce an auxiliary structure, which enables us to deal with quasiregular coherent configurations and generalized free products in a uniform way. To do this, we define the intersection of two sections $S=U/L$ and  $S'=U'/L'$ of a group $G$ by the formula
$$
S\cap S'=\css
(U\cap U')/LL'  &\text{if $LL'\le U\cap U'$},\\
\varnothing     &\text{otherwise}.\\
\ecss
$$
Note that the intersection is nonempty if $U=U'=G$ or $L=L'=1$. We say that~$S$ is a subsection of $S'$ if $L'\le L$ and $U\le U'$; in this case, we use notation $S\prec S'$. Clearly, the nonempty intersection $S\cap S'$ is a subsection of both $S$ and $S'$.\medskip

Let $I$ be a finite set, and let
\nmrt
\tm{F1} $\fG$ be a family of finite groups $G_i$ with identity~$e_i$, 
\tm{F2} $\fS$ be a family of sections $S_{ij}$ of $G_i$, 
where $S_{ii}=G_i/e_i$,
\tm{F3} $\fF$ be a family of isomorphisms 
$f_{ij}\in\iso(S_{ij},S_{ji})$, 
where $f_{ii}=f_{ij}f_{ji}=\id$,
\enmrt
where the indices $i$ and $j$ run over $I$. We assume that 
\qtnl{221216y}
S_{ij}\cap S_{ik}\ne\varnothing,\qquad i,j,k\in I.
\eqtn
This enables us to define the section $S_{ijk}=S_{ij}\cap S_{ik}$ of the group $G_i$; it is easily seen that $S_{ijk}=S_{ikj}$ for all $i,j,k$, and  $S_{ijk}=S_{ij}$ whenever $k=i$ or $k=j$. The isomorphism $f_{ik}$ takes the subsections of $S_{ik}$ to those of $S_{ki}$, and also preserves the partial order $\prec$. Therefore, $f_{ik}(S_{ijk})\prec f_{ik}(S_{ik})=S_{ki}$. We assume that 
\qtnl{221216z}
f_{ik}(S_{ijk})=S_{kij},\qquad i,j,k\in I.
\eqtn
Thus, the isomorphism $f_{ik}$ induces an isomorphism from $S_{ijk}$ onto $S_{kij}$; we denote it by $f_{ijk}$.

\dfntn
{\rm The triple $(\fG,\fS,\fF)$ defined by (F1), (F2), and (F3) is called a {\it system of linked sections} (based on $\fG$) if, in addition to conditions~\eqref{221216y} and~\eqref{221216z}, }
\qtnl{050817a}
f_{jki}f_{kij}f_{ijk}=\id,\qquad i,j,k\in I.
\eqtn
\edfntn	

The two special cases we are interested in are the systems of linked {\it quotients} and linked {\it subgroups}. In the former case, $S_{ij}=G_i/L_{ij}$ for all~$i,j$; in particular, $L_{ij}$ is a normal subgroup of~$G_i$. In the latter case $S_{ij}=U_{ij}/e_i$
for all~$i,j$; in this case, $S_{ij}$ is identified with a subgroup of $G_{ij}$. As we will see in Section~\ref{040816a}, a systems of linked quotients is nothing else than a quasiregular coherent configuration.\medskip

Let $(\fG,\fS,\fF)$ be a system of linked subgroups. In this case, $S_{ij}$ is identified with a subgroup $U_{ij}$ of the group $G_i$, and also $U_i:=U_{ii}$ is equal to~$G_i$, $i,j\in I$. Assume that all the groups $U_i$ are abelian. Then the direct product
$$
U=\prod_{i\in I}U_i
$$
is also abelian, and coincides with the free product of the groups $U_i$. It contains a subgroup
$$
U_0=\grp{x^{-1}f_{ij}(x):\ x\in U_{ij},\ i,j\in I},
$$
where $f_{ij}:U_{ij}\to U_{ji}$ is the isomorphism in (F3). Let
$$
\pi:U\to U/U_0
$$ 
be the canonical epimorphism. In this notation, the group $\pi(U)$ is called  the {\it generalized free product of the groups $U_i$  with amalgamated subgroups $U_{ij}$} if 
\qtnl{311216e}
\pi(U_i)\cong U_i\qaq\pi(U_i)\cap\pi(U_j)=\pi(U_{ij}),\qquad i,j\in I,
\eqtn
see paper \cite[p.~672]{N50}. It was remarked there that formulas~\eqref{050817a} form a necessary condition for existence of the generalized free product of the groups $A_i$  with amalgamated subgroups $A_{ij}$. In this sense, the system $(\fG,\fS,\fF)$ of linked subgroups will also be called an {\it amalgam configuration} (based on~$\fG$), and we say that it admits an amalgam if the generalized free product of the groups $U_i$  with amalgamated subgroups $U_{ij}$ does exists, i.e., condition~\eqref{311216e} is satisfied.\medskip

{\bf Duality for systems of linked sections.} The aim of this part is to define a dual of a system of linked sections based on a family of abelian groups. To this end, we recall  some standard facts from the duality theory for abelian groups and fix notations.\medskip

For an abelian group $G$, we denote by $\wh G$ its dual group, i.e. the group of complex-valued characters of $G$ with respect to componentwise multiplication. There is a canonical mapping from  the subgroup lattice of~$G$ to that of $\wh G$, taking a group $H\le G$ to the group
$$
H^\bot=\{\chi\in\wh G:\ H\subseteq  \ker(\chi)\},
$$ 
which is isomorphic to the dual group of $G/H$. Identifying the group $G$ with the dual group of $\wh G$, we have $(H^\bot)^\bot=H$. The mapping $\bot$ is an anti-isomorphism of the subgroup lattices. Namely, for all subgroups $L$ and $U$ of the group $G$,
\qtnl{050817y}
L\le U\quad\Leftrightarrow\quad U^\bot\le L^\bot
\eqtn  
and
\qtnl{271216a}
(LU)^\bot=L^\bot\cap U^\bot\qaq (L\cap U)^\bot=L^\bot U^\bot.
\eqtn
We extend the mapping $\bot$ to the sections of $G$ by setting
$(U/L)^\bot=L^\bot/U^\bot$. Clearly, this group is canonically isomorphic
to the dual group of $U/L$. It is a routine to check that for any two sections $S$ and $T$ of the group $G$, we have $S\cap T\ne\varnothing$ if and only if $S^\bot\cap T^\bot\ne\varnothing$, and in this case
\qtnl{311216t}
(S\cap T)^\bot=S^\bot\cap T^\bot.
\eqtn

Now let $U$ and $L$ be arbitrary subgroups of the group $G$, and let $\pi:U\to L$ a group homomorphism. Then one can define the mapping $\pi^\bot:L^\bot\to U^\bot$ that takes the character $\chi\in L^\bot$ to the character $\chi^{\pi^\bot}\in U^\bot$ defined by the formula
\qtnl{160817a}
\chi^{\pi^\bot}(u)=\chi(\pi(u)),\quad u\in U.
\eqtn
One can see that the mapping $\pi\mapsto\pi^\bot$ is a group isomorphism from $\Hom(U,L)$ onto $\Hom(L^\bot,U^\bot)$ that preserves isomorphisms. It is routine to define the isomorphism $\pi^\bot:T^\bot\to S^\bot$ for any isomorphism $\pi:S\to T$, where $S$ and $T$ are the sections of~$G$.\medskip

Let $(\fG,\fS,\fF)$ be a system of linked sections such that the family $\fG$ consists of abelian groups. Denote by $A$ the group dual to the group $G=\prod_{i\in I}G_i$. Then using formulas~\eqref{271216a}, one can see that $A$ is the direct product of the groups
$$
A_i=G_{i'}^\bot,\qquad i\in I,
$$
where $G_{i'}$ is the direct product of the $G_j$ with $j\ne i$. Assume  that $S_{ij}=U_{ij}/L_{ij}$, where $U_{ij}$ and $L_{ij}$ are subgroups of~$G_i\le G$, and also $L_{ij}\le U_{ij}$. Then by formulas~\eqref{271216a}, the group
 \qtnl{311216u}
A_{ij}:=((G_{i'}U_{ij})/(G_{i'}L_{ij}))^\bot=(A_i\cap L_{ij}^\bot)/(A_i\cap U_{ij}^\bot)
\eqtn
is a section of the group $A_i$ for all $i,j\in I$. By condition~\eqref{221216y}, the intersection $S_{ijk}=S_{ij}\cap S_{ik}$ is not empty. Therefore $S_{ijk}=U_{ijk}/L_{ijk}$, where $U_{ijk}=U_{ij}\cap U_{ik}$ and $L_{ijk}=L_{ij}L_{ik}$. Thus the set
$$
A_{ijk}:=A^{}_{ij}\cap A^{}_{ik}=((G_{i'}U_{ijk})/(G_{i'}L_{ijk}))^\bot
$$
is also not empty for all~$i,j,k$. By condition~\eqref{221216z}, the  
isomorphism $f_{ij}$ takes $U_{ijk}$ and $L_{ijk}$ to, respectively, $U_{kij}$ and $L_{kij}$. It induces an isomorphism 
$$
\wt f_{ik}:(G_{i'}U_{ijk})/(G_{i'}L_{ijk})\to(G_{k'}U_{kij})/(G_{k'}L_{kij}).
$$
Set $g^{}_{ik}=\wt f_{ki}^\bot$ for all $i,k\in I$. Then $g_{ik}(A_{ijk})=A_{kij}$ and $g_{ijk}g_{kij}g_{jki}=\id$ for all $i,j,k\in I$. Thus, the triple 
\qtnl{241216t}
(\fG,\fS,\fF)^\bot:=(\fG^\bot,\fS^\bot,\fF^\bot)
\eqtn
where $\fG^\bot=\{A_i\}_{i\in I}$, $\fS^\bot=\{A_{ij}\}_{i,j\in I}$, and $\fF^\bot=\{g_{ij}\}_{i,j\in I}$, satisfies conditions~\eqref{221216y}, \eqref{221216z} and \eqref{050817a}. This proves the first part of the following statement (the second part follows by duality).

\lmml{241216r}
The triple~\eqref{241216t} is a system of linked sections. Moreover, 
$$
((\fG,\fS,\fF)^\bot)^\bot=(\fG,\fS,\fF).\footnote{Here, we use canonical isomorphisms between $A_i^\bot$ and $G^{}_i$, and between $A_{ij}^\bot$ and $G_{ij}^{}$.}
$$
\elmm

Note that if $(\fG,\fS,\fF)$ is a system of linked quotients (respectively, subgroups), then $(\fG,\fS,\fF)^\bot$ is a system of linked subgroups (respectively, quotients). This fact together with Lemma~\ref{241216r} immediately imply the following statement.

\crllrl{241216q}
For a fixed family $\fG$ of finite abelian groups, the mapping 
$$
(\fG,\fS,\fF) \to (\fG,\fS,\fF)^\bot
$$
induces a one-to-one correspondence between the systems of linked quotients based on~$\fG$ and the systems of linked subgroups based on~$\fG^\bot$.
\ecrllr

\section{Quasiregular coherent configurations as systems of linked quotients}\label{040816a}
Throughout this section, $\cX=(\Omega,S)$ is a quasiregular coherent configuration with fibers $\Omega_1,\ldots,\Omega_m$, $m\ge 1$, and $I=\{1,\ldots,m\}$. For each $i\in I$, the $i$th  homogeneous component $\cX_i$ of $\cX$ is a (regular) scheme the basis relations of which form a group~$G_i$, and we say that $\cX$ is of type $\fG=\{G_i\}_{i\in I}$. Note that any quasiregular coherent configuration is of type $\fG$ for suitable family of groups $G_i\in\fG$. The aim of this section is to prove that a quasiregular coherent configuration of type $\fG$ is nothing else than a system of linked quotients based on~$\fG$.\medskip

Let $i,j\in I$. In what follows, we set $S_{ij}:=S_{\Omega_i,\Omega_j}$; in particular, $S_{ii}=G_i$. Note that here, $S_{ij}$ is a subset of $S$ and not a section of a group as in Section~3. Since any relations from $G_i$ and $G_j$ are thin, we have 
$$
s_i^*S_{ij}=S_{ij}=S_{ij}s^{}_j,\qquad s_i\in G_i,\ s_j\in G_j.
$$
This implies that the mappings $\pi_{ij}:G_i\to\sym(S_{ij})$
and $\rho_{ij}:G_j\to\sym(S_{ij})$ given by the formulas
$$
x^{\pi_{ij}(s)}=s^*x,\ x\in S_{ij},\qaq 
x^{\rho_{ij}(s)}=xs,\ x\in S_{ij},
$$
are group homomorphisms. 

\lmml{tran}
In the above notation, for all $i,j\in I$,
\nmrt
\tm{1} $\im(\pi_{ij})$ and $\im(\rho_{ij})$ are regular groups centralizing each other,
\tm{2} $\ker(\pi_{ij})=xx^*$ for all $x\in S_{ij}$, and $\ker(\rho_{ij})=x^*x$ for all $x\in S_{ji}$.
\enmrt
\elmm
\proof Let $x,y\in S_{ij}$. Then any basis relations $s_i\in xy^*$ and $s_j\in x^*y$ belong to $G_i$ and $G_j$, respectively. Therefore they are thin and hence
$$
s_i^*x=y=xs^{}_j. 
$$
This proves that $\im(\pi_{ij})$ and $\im(\rho_{ij})$ are transitive groups. Obviously each of them centralizes the other.  Next assume that $s_i^*x=x$ for some $x\in S_{ij}$. Since the group $\im(\rho_j)$ is transitive, any $y\in S_{ij}$ is of the form $xs_j$ for a suitable $s_j$. Therefore,
$$
s_i^*y=s_i^*xs_j^{}=xs_j^{}=y.
$$
Thus, the group $\im(\pi_{ij})$ is regular. Similarly, one can prove that the group $\im(\rho_{ij})$ is also regular. This proves statement~(1). Now by the regularity, $s\in \ker(\pi_{ij})$ if and only if $s^*x=x$ for some (and hence for all) $x\in S_{ij}$ if and only if $s\in xx^*$. Similarly, $s\in \ker(\rho_{ij})$ if and only if $s\in x^*x$ for some (and hence for all) $x\in S_{ij}$.\eprf\medskip

In what follows, we set 
$$
G_{ij}=\ker(\pi_{ij})=\ker(\rho_{ji}),\qquad i,j\in I.
$$
From statement~(1) of Lemma~\ref{tran}, it follows that the groups $G_i/G_{ij}$ and $G_j/G_{ji}$ are the left and right representations of the same group. In particular, these groups are isomorphic. To write an explicit isomorphism, choose arbitrarily points $\alpha_i\in \Omega_i$, $i\in I$, and set
$$
s_{ij}:=r(\alpha_i,\alpha_j),\quad i,j\in I.
$$ 

\lmml{100817a}
For all $i,j\in I$, the mapping
$$
f_{ij}:G_i/G_{ij}\to G_j/G_{ji},\ G_{ij}s_i\mapsto s_{ij}^*(G^{}_{ij}s_i^{})s^{}_{ij}
$$
is a group isomorphism, and $(f_{ij})^{-1}=f_{ji}^{}$.
\elmm
\proof To prove that $f_{ij}$ is well-defined, denote $s_{ij}$ by $x$. Then $s_ix=xs_j$ for some $s_j\in G_j$ and $G_{ij}=xx^*$ by statements~(1) and~(2) of Lemma~\ref{tran}, respectively. Thus, by Lemma~\ref{sss}, we have
\qtnl{110817a}
x^*( G_{ij}s_i)x=x^* (xx^*)s_ix=x^*s_ix=x^*xs_j=G_{ji}s_j\in G_j/G_{ij},
\eqtn
which shows that $f_{ij}$ is well-defined. It is a homomorphism, because for any elements $s_i,t_i\in G_i$,
$$
f_{ij}(G_{ij}s_i\,G_{ij}t_i)=f_{ij}(G_{ij}s_it_i)=x^*(s_it_i)x=
(x^*s_ix)(x^*t_i x)=f_{ij}(G_{ij}s_i)f_{ij}(G_{ij}t_i).
$$
Assume that $f_{ij}(G_{ij}s_i)=G_{ji}$. Then formula~\eqref{110817a}
shows that $s_j\in G_{ji}$, and hence $s_ix=xs_j\in xx^*x=x$.
It follows that $s_i^*\in G_{ij}$, and therefore $s_i\in G_{ij}$, which means that the homomorphism $f_{ij}$ is one-to-one. Formula~\eqref{110817a} also shows that it is an epimorphism, because for every $s_j\in G_j$, there exists $s_i\in G_i$ such that $s_ix=xs_j$ (statement~(1) of Lemma~\ref{tran}). Thus, $f_{ij}$ is an isomorphism. Since also
$$
f_{ij}f_{ji}(G_{ji}s_i)=f_{ij}(xs_ix^*)=x^*xs_ix^*x=G_{ji}s_i,
$$
we conclude that $f_{ij}f_{ji}$ is the identity map.\eprf\medskip

At this point, starting with the quasiregular coherent configuration $\cX$ with fibers indexed by the set $I$, we have constructed a triple 
\qtnl{140817a}
\cT=\cT(\cX)=(\{G_i\}_{i\in I}, \{G_i/G_{ij}\}_{i,j\in I}, \{f_{ij}\}_{i,j\in I}).
\eqtn

\lmml{110817b}
The triple $\cT$ is a system of linked quotients.
\elmm
\proof To verify condition~\eqref{221216z}, let $i,j,k\in I$. From the definition of $s_{ij}$, it follows that $s_{jk}\in s_{ji}s_{ik}$. This implies that
$$
G_{jk}=s_{jk}^{}s_{jk}^*\subseteq (s_{ji}s_{ik})(s_{ji}s_{ik})^*=s_{ij}^* s_{ik}^{}s_{ik}^*s_{ij}^{}=f_{ij}(G_{ik}),
$$
and similarly, $G_{ik}\subseteq f_{ji}(G_{jk})$. Therefore,
$f_{ij}(G_{ik})\subseteq f_{ij}f_{ji}(G_{jk})=G_{ji}G_{jk}$.
Thus, 
$$
G_{ji}G_{jk}=f_{ij}(G_{ik}G_{ij}).
$$
This enables us to define the induced (by $f_{ij}$) isomorphism
$$
f_{ijk}:G_i/G_{ij}G_{ik}\to G_j/G_{ji}G_{jk}.
$$
It remains to verify condition~\eqref{050817a}, or, equivalently, that the composite $f_{kij}f_{jki}f_{ijk}$ is the identity map. To this end, let $x\in G_{ii}$. Then
\begin{eqnarray*}
f_{kij}f_{jki}f_{ijk}(G_{ij}G_{ik}x) &=& f_{kij}f_{jki}(s_{ji}xs_{ij}G_{jk})\\
& =&f_{kij}(s_{kj}s_{ji}xs_{ij}s_{jk})\\ 
& =&s_{ik}s_{kj}s_{ji}xs_{ij}s_{jk}s_{ki}.
\end{eqnarray*}
Note that $1_{\Omega_i}\in s_{ij}s_{jk}s_{ki}\cap s_{ik}s_{kj}s_{ji}$ by the definition of the relations $s_{ij}$. Therefore the left-hand side of the above equality equals the coset of $G_{ij}G_{ik}$ that contains~$x$,
and the right-hand side must be one coset, i.e., the coset $(G_{ij}G_{ik}x)$, as required.\eprf\medskip

Let us show that every system of linked quotients based on a family $\fG$ of groups $G_i$, $i\in I$, is of the form $\cT(\cX)$ for some quasiregular coherent configuration $\cX$ of type~$\fG$. To this end, assume that $I$ is a finite set and
$$
(\fG,\fS,\fF)=(\{G_i\}_{i\in I}, \{G_{ij}\}_{i,j\in I},\{f_{ij}\}_{i,j\in I})
$$ 
is a system of linked quotients. Denote by $\Omega$ the disjoint union of all the~$G_i$. For each $i\in I$ and each $x\in G_i$, we define a binary relation $s_x$ on the set $\Omega_i=G_i$ as the graph of the left multiplication by~$x$,
$$
s_x=\{(\alpha,\beta)\in \Omega_i\times \Omega_j:\ \beta=x\alpha\}.
$$

\lmml{120817f}
Given indices $i,j\in I$, the set $\Omega_i\times\Omega_i$ is partitioned into the relations $s_x\cdot s_{ij}$, $x\in G_i$, where
$$
s_{ij}=\bigcup_{g_i\in G_i}\{g_i\}\times f_{ij}(G_{ij}g_i).
$$
\elmm
\proof  We have to verify that if $x,y\in G_i$, then the relations $s_x\cdot s_{ij}$ and $s_y\cdot s_{ij}$ are disjoint or coincide. To this end, we note that
$$
s_x\cdot s_{ij}=\bigcup_{\gamma\in G_i}\{x\gamma\}\times f_{ij}(G_{ij}\gamma)=\bigcup_{\delta\in G_i}\{\delta\}\times f_{ij}(G_{ij}x^{-1}\delta).
$$
Now if the relations $s_x\cdot s_{ij}$ and $s_y\cdot s_{ij}$ are intersected, say in a pair $(\alpha,\beta)$, then this implies that $\beta$ belongs to the intersection of two right cosets  $f_{ij}(G_{ij}x^{-1}\alpha)$ and $f_{ij}(G_{ij}(y^{-1}\alpha)$ of the same group $G_{ji}\le G_j$. It follows that $G_{ij}x^{-1}=G_{ij}y^{-1}$ and hence $x=y z$ for some $z\in G_{ij}$. Since $z\cdot s_{ij}=s_{ij}$, we conclude that 
$$
s_x\cdot s_{ij}= s_{yz}s_{ij}=s_y\cdot(s_z s_{ij})=s_y\cdot s_{ij}
$$
as required.\eprf\medskip

Denote by $S$ the union of the sets $S_{ij}=\{s_xs_{ij}\mid x\in G_i\}$ over all $i,j\in I$.

\lmml{120817a}
The pair $\cX=(\Omega, S)$ is a quasiregular coherent configuration such that $\cT(\cX)=(\fG,\fS,\fF)$.
\elmm
\proof We have already proved that $S$ is a partition of $\Omega\times \Omega$. It is easily verified that  $1_\Omega$ is a disjoint union of $s_{e_i}$ with $i\in I$, where $e_i$ is the identity of $G_i$. Next, let $i,j\in I$. By the definition of $s_{ij}$, we have
$$
s_{ij}^*=\bigcup_{\gamma\in G_i} f_{ij}(G_{ij}\gamma)\times \{\gamma\}.
$$
Take $(\alpha,\beta)\in s_{ij}^*$. Then $\beta\in f_{ij}(G_{ij}\alpha)$, which implies $f_{ji}(G_{ji}\beta)=G_{ij}\alpha$, because $(f_{ij})^{-1}=f_{ji}$. Thus, 
$$
(\beta,\alpha)\in \{\beta\}\times f_{ji}(G_{ji}\beta)\subseteq s_{ji}
$$
and, hence $s_{ij}^* \subseteq s_{ji}$. Since $i$ and $j$ were taken arbitrarily, it follows that $s_{ji}^*\subseteq s_{ij}$ and hence $s_{ji}\subseteq s_{ij}^*$. Thus, $s_{ij}^*=s_{ji}$ and so $S^*=S$.\medskip

We claim that $s_{ij}s_{jk}$ equals the union of all sets $s_xs_{ik}$
with $x\in G_{ij}$. To this end, let $(\alpha,\gamma)\in s_{ij}s_{jk}$.
Then there exists $\beta\in G_j$ such that $(\alpha,\beta)\in s_{ij}$ and $(\beta,\gamma)\in s_{jk}$. In particular, we have $\beta\in f_{ij}(G_{ij}\alpha)$ and $\gamma\in f_{jk}(G_{jk}\beta)$. Since $G_{jk}G_{ji}=f_{ij}(G_{ij}G_{ik})$ and $f_{jki}f_{ijk}f_{kij}=\id$, we conclude that
\begin{eqnarray*}
\gamma\in f_{jk}(G_{jk}f_{ij}(G_{ij}\alpha))& =&f_{jk}(f_{ij}(G_{ij}G_{ik})f_{ij}(G_{ij}\alpha))\\
 & =& f_{jk}f_{ij}(G_{ij}G_{ik}\alpha)\\
 & =&f_{ki}^{-1}(G_{ij}G_{ik}\alpha)\\
 &=&f_{ik}(G_{ij}G_{ik}\alpha).
\end{eqnarray*}
This implies that $(\alpha x,\gamma)\in s_{ik}$ for some $x\in G_{ij}$, and hence,
$$
(\alpha,\gamma)\in  s_xs_{ik}\subseteq \bigcup_{y\in G_{ij}}s_ys_{ik}.
$$
Thus, $s_{ij}s_{jk}\subseteq\bigcup_{x\in G_{ij}}s_xs_{ik}$. Conversely, if $x\in G_{ij}$, then 
$$
s_xs_{ik}\subseteq s_{ij}s_{ji}s_{ik}\subseteq s_{ij}\bigcup_{y\in G_{ji}}s_ys_{jk}=s_{ij}s_{jk},
$$
which completes the proof of the claim.\medskip

It is easily verified that for each $x\in G_i$, there exists $y\in G_j$ such that $s_xs_{ij}=s_{ji}s_y$. Furthermore, from the claim it follows that $s_{ij}\cdot s_{jk}$ is invariant with respect to the left multiplication by $s_z$ with $z\in G_{ij}G_{ik}$. Therefore the product of the adjacency matrices of $s_xs_{ij}$ and $s_{jk}s_y$ is a scalar multiple of the adjacency matrix of $s_{ij}\cdot s_{jk}$. This proves the existence of the intersection numbers for the partition~$S$. Thus, $\cX$ is a coherent configuration. The second statement of the lemma is straightforward.\eprf

\crllrl{220817a}
Let $\cX$ be a quasiregular coherent configuration of type~$\fG=\{G_i\}_{i\in I}$. Then there exist groups $G_{ij}\trianglelefteq G_i$, $i,j\in I$, such that any basis relation $s\in S_{ij}$ is of the form 
$$
s=\bigcup_{y\in G_i}\{x_sy\}\times f_{ij}(G_{ij}y).
$$
for a suitable $x_s\in G_i$.
\ecrllr

From Lemmas~\ref{110817b} and~\ref{120817a}, we immediately get the main result of this section.

\thrml{120817b}
Let $\fG$ be a family of finite groups. Then the mapping $\cX\to\cT(\cX)$ gives a one-to-one correspondence between the quasiregular coherent configurations of type~$\fG$ and the systems of linked quotients based on $\fG$.
\ethrm

As an example, let us consider a special case when all the groups of the family~$\fG$ are simple. Then for all $i,j\in I$, the group $G_{ij}$ being a normal subgroup of $G_i$ is either $G_i$ itself or~$1$. Moreover, if $G_{ij}=G_{jk}=1$ for some $k\in I$, then obviously $G_{ik}=1$. This implies that the relation on $I$, defined by
$$
i\sim j\quad\Leftrightarrow\quad G_{ij}=1,
$$
is an equivalence relation, and if $i\sim j$, then the groups $G_i$ and $G_j$ are isomorphic and $|S_{ij}|=|S_{ji}|=|G_i|=|G_j|$. Therefore, if $J$ is a class of $\sim$ and $\Omega_J$ is the union of all $\Omega_i$ with $i\in J$, then the coherent configuration $\cX^{(J)}=(\Omega_J,S_J)$ is semiregular, where $S_J$ is the union of all $S_{ij}$ with $i,j\in J$.
On the other hand, 
$$
i\not\sim j\quad\Rightarrow\quad S_{ij}=\{\Omega_i\times\Omega_j\}.
$$
Therefore, the coherent configuration $\cX$ is the direct sum of the coherent configurations $\cX^{(J)}$, where $J$ runs over the classes of the equivalence relation~$\sim$. According to \cite[Theorem~3.3]{EP09} each of these coherent configuration is schurian and separable. Since the direct sum of schurian (respectively, separable)  coherent configurations is also schurian (respectively, separable) \cite[Theorem~7.5]{EP09}, this proves the following statement.

\thrm
Let $\fG$ be a family of finite simple groups. Then any quasiregular coherent configuration $\cX$ of type~$\fG$ is the direct sum of semiregular
coherent configurations. In particular, $\cX$ is schurian and separable.
\ethrm

\section{Proof of Theorem~\ref{241216g}}\label{220817c}
Let $\fG=\{G_i\}_{i\in I}$, where $G_i$ is an abelian group and $I$ is a finite set. Then by Theorem \ref{120817b} and Corollary \ref{241216q}, the mapping
$$
\cX\mapsto\cT(\cX)^\bot
$$
defines a one-to-one correspondence between the quasiregular coherent configurations of type $\fG$ and the amalgam configurations based on~$\fG$.
This proves the first part of the theorem.\medskip

To prove the second part, let $\cX=(\Omega,S)$ be a quasiregular coherent configuration of type~$\fG$, and let $\cT=\cT(\cX)$ be the associated system of linked quotients defined by~\eqref{140817a}. Define the group
$$
G=\prod_{i\in I}G_i
$$
considered as a permutation group on $\Omega$ such that $\alpha s=\{\alpha^s\}$ for all $\alpha\in\Omega_i$, $s\in G_i$, and $i\in I$. 
In what follows, for each $i\in I$ we denote by $G_{i'}$ the direct product of the groups $G_j$ with $j\ne i$. 

\lmml{140811a}
The coherent configuration $\cX$ is schurian if and only if there exists a group $H\le G$ such that for all $i,j\in I$,
\qtnl{291216a}
H^{\Omega_i\cup\Omega_j}=\{(s_i,s_j)\in G_i\times G_j:\ f_{ij}(G_{ij}s_i)=G_{ji}s_j\}.
\eqtn
\elmm
\proof From Lemma~\ref{120817f}, it follows that 
the condition given in the lemma is equivalent to the fact that $\orb(H,\Omega_i\times\Omega_j)=S_{ij}$. Since the latter is true for all $i,j\in I$ if and only if $\cX$ is schurian, we are done.\eprf\medskip

To complete the proof, we verify that if $\cX$ is schurian, then the amalgam configuration $\cT^\bot$ admits an amalgam (the converse statement can be proved in a similar way). To this end, assume that $\cX$ is schurian. Then condition~\eqref{291216a} is satisfied for some group $H\le\sym(\Omega)$; we  take $H=\aut(\cX)$. In what follows, assume that  $\cT^\bot$ is the amalgam configuration as in formula~\eqref{241216t}.\medskip 

Let $i,j\in I$. If $i=j$, then condition~\eqref{291216a} implies that $H^{\Omega_i}=G_i$, or equivalently, that $HG_{i'}=G$. Passing to the dual subgroups, we get from this and the first equality in~\eqref{271216a} that $A_i\cap H^\bot=1$ for all $i\in I$. Therefore if $\pi:A\to A/H^\bot$ is the canonical homomorphism, then
\qtnl{311216x}
\pi(A_i)=A_iH^\bot,\qquad i\in I.
\eqtn
Now assume that $i$ and $j$ are not necessarily equal. From condition~\eqref{291216a}, it follows that the group $H^{\Omega_{i}\cup\Omega_j}$ contains a subgroup $G_{ij}\times G_{ji}$. It follows that
$$
G_{ij}G_{i'}\cap H=((G_{ij}\times G_{ji})(G_{i'}\cap G_{j'}))\cap H=G_{ji}G_{j'}\cap H.
$$
Note that in formula~\eqref{311216u} written for the amalgam configuration $\cT^\bot$, we should have $U_{ij}=G_i$ and $L_{ij}=G_{ij}$. Therefore, $A_{ij}=A_i\cap G_{ij}^\bot$ and $A_{ji}=A_i\cap G_{ji}^\bot$. It follows that passing to the dual subgroups in the left- and right-hand sides of the above equality, we have
\qtnl{311216y}
\pi(A_{ij})=A_{ij}H^\bot=(A_i\cap G_{ij}^\bot)H^\bot=(A_j\cap G_{ji}^\bot)H^\bot=A_{ji}H^\bot=\pi(A_{ji}).
\eqtn
for all $i,j\in I$. Formulas~\eqref{311216x} and~\eqref{311216y}  imply condition~\eqref{311216e} in the definition of generalized free product of the groups $A_i$  with amalgamated subgroups $A_{ij}$. Thus, it remains to verify that the kernel $H^\bot$ of the homomorphism $\pi$ is generated by the elements $x^{-1}g_{ij}(x)$, where $x\in A_{ij}$ and $i,j\in I$. To this end, we note that that formula~\eqref{291216a} shows that
\qtnl{311216i}
\aut(\cX)=\bigcap_{i,j\in I}M_{ij},
\eqtn
where $M_{ij}=\{(g_i)_{i\in I}\in G:\ f_{ij}(G_{ij}g_i)=G_{ji}g_j\}$. 
However, using the definitions of the operator $\bot$ and isomorphism $g_{ij}$, one can easily check that the group $M_{ij}^\bot$ consists of all elements $a^{-1}g_{ij}(a)$
with $a\in A_i$. Taking into account that
$$
\aut(\cX)^\bot=\langle M_{ij}^\bot:\ i,j\in I\rangle
$$
(see formula~\eqref{311216i}), we arrive to the required statement on the group $H^\bot=\aut(\cX)^\bot$.

\section{Proof of Theorem~\ref{130617a}}\label{220817y}

In this section, we deduce Theorem~\ref{130617a} from a quite general result  establishing schurity and separability for a large class of coherent configurations properly containing the class studied in~\cite{Zies}. Namely,
let $\cX=(\Omega,S)$ be a {\it coset}\,\footnote{The name comes from \cite{EP16}, where the coherent configurations admitted a regular abelian automorphism group and satisfying condition~\eqref{180817a} were studied under that name.} coherent configuration, i.e., the following condition holds:
\qtnl{180817a}
ss^*s=s\quad\text{for all}\ \,s\in S.
\eqtn
From Lemma~\ref{sss}, it follows that every quasiregular coherent configuration is a coset one. We say that a coset  coherent configuration $\cX$ satisfies a {\it distributivity condition} if 
$$
W\cap UV=(W\cap U)(W\cap V)\quad\text{for all}\ \,U,V,W\in\cL(\Delta),\ \Delta\in F,
$$
where $F=F(\cX)$ and $\cL(\Delta)$ is the set of all intersections of closed subsets $ss^*$ with $s\in S_{\Delta,\Gamma}$ and $\Gamma\in F$.  This condition obviously holds for any quasiregular coherent configuration whose homogeneous components correspond the groups with distributive lattice of normal subgroups.

\lmml{faith}
Let $\cX$ be a coset coherent configuration satisfying the distributivity condition, and let $\varphi$ be an algebraic isomorphism from $\cX$ to a coherent configuration~$\cX'$. Then
\nmrt
\tm{1} $\cX'$ is a coset coherent configuration satisfying the distributivity condition, 
\tm{2} each $\varphi$-faithful map is $\varphi$-extendable.
\enmrt
\elmm
\proof Statement~(1) is obvious. To prove statement~(2), let us verify that a $\varphi$-faithful map $f$ is $\varphi$-extendable to a point~$\gamma\in\Omega$. Without loss of generality we assume that the set  $\Lambda=\dom(f)$ contains at least three distinct elements two of which, say $\alpha$ and $\beta$. Set 
$$
\Lambda_\alpha=\Lambda\setminus\{\beta\}\qaq
\Lambda_\beta=\Lambda\setminus\{\alpha\},
$$
in particular, $\alpha\in\Lambda_\alpha$ and $\beta\in\Lambda_\beta$. By induction on $|\dom(f)|$, we may assume that the restrictions $f_\alpha$ and $f_\beta$ of the map $f$ to, respectively, $\Lambda_\alpha$ and $\Lambda_\beta$ are $\varphi$-faithful maps. In view of formula~\eqref{100617u}, this implies that
\qtnl{120617a}
\bigcap_{\lambda\in\dom(f_\alpha)}\lambda^{f_\alpha}(r_{\lambda\gamma})^\varphi\ne\varnothing\qaq
\bigcap_{\lambda\in\dom(f_\beta)}\lambda^{f_\beta}(r_{\lambda\gamma})^\varphi\ne\varnothing,
\eqtn
where $r_{\lambda\gamma}=r(\lambda,\gamma)$. Take arbitrary elements $\alpha'$ and~$\beta'$ belonging to the first and second intersection, respectively. Then for all $\lambda\in\Lambda_\alpha\cap\Lambda_\beta$, we have
$$
r'(\alpha',\beta')\in (r_{\gamma\lambda})^\varphi (r_{\lambda\gamma})^\varphi=(r_{\gamma\lambda}r_{\lambda\gamma})^\varphi\subseteq W^\varphi=:W',
$$
where $W$ is the intersection of all sets $r_{\gamma\lambda}r_{\lambda\gamma}$ with $\lambda\in\Lambda_\alpha\cap\Lambda_\beta$. Note that since $\cX$ is a coset coherent configuration, $W'$ is a closed subset of $\cX'_{\Gamma'}$, where $\Gamma$ is the fiber containing~$\gamma$ and $\Gamma'=\Gamma^\varphi$. Next, 
$$
r'(\alpha^f,\beta^f)=(r_{\alpha\beta})^\varphi\in(r_{\alpha\gamma}r_{\gamma\beta})^\varphi=(r_{\alpha\gamma})^\varphi (r_{\gamma\beta})^\varphi.
$$
Consequently, there exists a point $\gamma'\in\Omega'$ such that 
$$
r'(\alpha^f,\gamma')=(r_{\alpha\gamma})^\varphi\qaq r'(\gamma',\beta^f)=(r_{\gamma\beta})^\varphi.
$$
Furthermore, from the definition of $\alpha'$ and formula~\eqref{120617a} for $\lambda=\alpha$, we have $r'(\alpha^f,\alpha')=(r_{\alpha\gamma})^\varphi$. Therefore,
$$
r'(\gamma',\alpha')\in r'(\gamma',\alpha^f)r'(\alpha^f,\alpha')=(r_{\gamma\alpha})^\varphi (r_{\alpha\gamma})^\varphi=(r_{\gamma\alpha}^{}r_{\gamma\alpha}^\ast)^\varphi=:U^\varphi,
$$
Similarly, 
$$
r'(\gamma',\beta')\in (r_{\gamma\beta}^{}r_{\gamma\beta}^\ast)^\varphi=:V^\varphi.
$$
Again, $U'=U^\varphi$ and $V'=V^\varphi$ are closed subsets of the coherent configuration $\cX'_{\Gamma'}$. Moreover,
$$
r'(\gamma',\alpha')\in U',\quad 
r'(\gamma',\beta')\in V',\quad
r'(\alpha',\beta')\in W'.
$$
Since also $W'\cap U'V'=(W'\cap U')(W'\cap V')$ (the distributivity condition), by the Zieschang lemma \cite[Lemma~1]{Zies} there exists a point $\ov\gamma\in\Omega'$ such that
$$
\ov\gamma\in \alpha'(W'\cap U')\,\cap\, \beta'(V'\cap W')\,\cap\, \gamma'(U'\cap V').
$$
To complete the proof, it suffices to verify that 
\qtnl{130617b}
\ov\gamma\in \bigcap_{\lambda\in\Lambda}\lambda^f(r_{\lambda\gamma})^\varphi.
\eqtn
Since 
$r_{\alpha\gamma}(U\cap V)$ is contained in $r_{\alpha\gamma}^{}r_{\alpha\gamma}^\ast r_{\alpha\gamma}^{}=r_{\alpha\gamma}^{}$,
we have
$$
r'(\alpha^f,\ov\gamma)\in r'(\alpha^f,)inp{\gamma'})r'(\gamma',\ov\gamma)
=(r_{\alpha\gamma})^\varphi(U\cap V)^\varphi=(r_{\alpha\gamma}(U\cap V))^\varphi=(r_{\alpha\gamma})^\varphi
$$

and hence, $\ov\gamma\in \alpha^f(r_{\alpha\gamma})^\varphi$.
Similarly, we obtain $\ov\gamma\in \beta^f(r_{\beta\gamma})^\varphi$.
Finally for any $\lambda\in \Lambda_\alpha\cap\Lambda_\beta$, we have
$$
r'(\lambda^f,\ov\gamma)\in r'(\lambda^f,\beta')r'(\beta',\ov\gamma)\subseteq
(r_{\lambda\gamma})^\varphi W^\varphi\subseteq(r_{\lambda\gamma}^{}r_{\lambda\gamma}^*r_{\lambda\gamma})^\varphi=(r_{\lambda\gamma})^\varphi
$$
which proves~\eqref{130617b}.\eprf\medskip

By the above remarks, Theorem~\ref{130617a} immediately follows from the theorem below, which is a straightforward consequence of Lemma~\ref{seps} and \ref{faith}.

\thrml{distributive}
Any coset coherent configuration satisfying the distributivity condition is  schurian and separable.
\ethrm

\section{Quasiregular coherent configurations with at most four fibers}\label{220817u}

Throughout this section, $\cX=(\Omega,S)$ is a quasiregular
coherent configuration with fibers $\Omega_i$ and $S_{ij}=S_{\Omega_i,\Omega_j}$, $i,j\in I$. First, we consider the case when $|I|\le 3$.\medskip

{\bf Proof of Theorem~\ref{120817d}.} Let $\cX'=(\Omega',S')$ be a coherent configuration, and let 
$$
\varphi:S\to S',\ s\to s'
$$  
be an algebraic isomorphism. It induces a bijection $\Omega_{i^{}}\mapsto\Omega'{}_{i'}$ between the fibers of $\cX$ and $\cX'$, and an algebraic isomorphism from $\cX^{}_{\Omega^{}_i}$ onto $\cX'_{\Omega'_{i'}}$, $i\in I$. In particular, the coherent configuration $\cX'$ is also quasiregular.\medskip

For each $i\in I$, we choose an arbitrary point $\alpha_i\in \Omega_i$, and set
$$
s_{ij}:=r(\alpha_i,\alpha_j),\qquad i,j\in I.
$$
Then obviously $s_{ij}\in s_{ik}s_{kj}$, and hence $s'_{ij}\in s'_{ik}s'_{kj}$ for all $i,j,k$. Since $|I|\le 3$, there exist points
$\alpha'_i\in \Omega'$, $i\in I$
such that 
$$
s'_{ij}=r(\alpha'_i,\alpha'_j),\qquad i,j\in I.
$$
Furthermore, taking into account that the the coherent configurations $\cX$ and $\cX'$ are quasiregular, we conclude that
\qtnl{170817q}
\Omega^{}=\{\alpha^{}_i s^{}:\ s^{}\in S^{}_{ii},\ i\in I\}\qaq 
\Omega'=\{\alpha'_i s':\ s'\in S'_{ii},\ i\in I\},
\eqtn
where the singletons $\alpha^{}_i s^{}$ and $\alpha'_i s'$ are identified with the corresponding points of $\Omega^{}_i$ and $\Omega'_i$, respectively.
It follows that the mapping
\qtnl{170817g}
f:\Omega\to \Omega',\ \alpha^{}_is\mapsto \alpha'_is^\varphi,
\eqtn
where $s\in S_{ii}$, is a bijection. We claim that it is  a combinatorial isomorphism from $\cX$ onto $\cX'$. Indeed, if $\alpha,\beta\in \Omega$, then in view of~\eqref{170817q}  there exist $s\in S_{ii}$ and $t\in S_{jj}$
such that $\alpha=\alpha_i s$ and $\beta=\alpha_jt$. We have
\begin{eqnarray*}
r(\alpha^f,\beta^f) & =&r(\alpha_i's^\varphi, \alpha_j't^\varphi)\\
& =&(s^\varphi)^* r(\alpha_i',\alpha_j')t^\varphi\\
& =&(s^\varphi_{})^* (s^{}_{ij})^\varphi t_{}^\varphi\\
& =&(s^* s_{ij}t)^\varphi=r(\alpha,\beta)^\varphi.
\end{eqnarray*}

Therefore the algebraic isomorphism $\varphi$ is induced by the combinatorial isomorphism~$f$. Thus, the coherent configurations $\cX$ is separable.\medskip

To prove the schurity of $\cX$, let $r(\alpha,\beta)=r(\alpha',\beta')$, where $\alpha',\beta'\in\Omega$. Assume first that $\alpha\in\Omega_i$ and $\beta\in\Omega_j$ with $i\ne j$. Then $(\alpha,\beta)^f=(\alpha',\beta')$, where $f$ is the isomorphism defined by \eqref{170817g} for $\varphi=\id_S$, $(\alpha_i,\alpha_j)=(\alpha,\beta)$,  and $(\alpha'_i,\alpha'_j)=(\alpha',\beta')$. It follows that for the group $K$ generated by these isomorphisms $f$ with all possible $(\alpha,\beta)$ and $(\alpha',\beta')$, we have
$$
\orb(K,\Omega_i\times\Omega_j)=S_{ij}
$$
for all distinct $i$ and $j$. It immediately implies that if $|I|\ne 1$, then $K$ acts transitively on each $\Omega_i$, and hence $\cX$ is schurian. Since the schurity of $\cX$ for $|I|=1$ is obvious, we are done.\eprf\medskip

{\bf Proof of Corollary~\ref{120817e}.} Let $\cX$ be a meta-thin scheme with thin residue of index at most~$3$. Denote by $\cY$ the thin residue extension of~$\cX$ in the sense of~\cite[Sec.~2.4]{MP}. Then $\cY$ is a quasiregular coherent configuration with at most~$3$ fibers (which are geometric cosets of the thin residue). By  Theorem~\ref{120817d}, it is schurian and separable. According to statement~(3) in~\cite[Theorem~3.3]{MP}, this shows that $\cX$ is schurian and separable.\eprf\medskip

In the rest of the section, we construct a huge class of non-schurian quasiregular coherent configurations with exactly four fibers. To this end, let $G=C_p\times C_p$, where $p$ is a prime, and let $\fG=\{G_i\}_{i\in I}$, where $I=\{1,2,3,4\}$ and $G_i=G$ for all $i\in I$. Choose a family $\cG$ of groups $G_{ij}\le G_i$, $i,j\in I$, such that
\qtnl{180817f}
|G_{ij}|=\css
p &\text{if $i\ne j$,}\\
1 &\text{if $i=j$.}\\
\ecss
\eqtn
Then the factor groups $G_i/G_{ij}$ and $G_j/G_{ji}$ are isomorphic for all~$i,j$. Therefore, one can choose a family $\fF$ of isomorphisms $f_{ij}:G_i/G_{ij}\to G_j/G_{ji}$, $i,j\in I$, such that
\qtnl{180817g}
f_{ii}=\id_G\qaq f_{ij}=(f_{ji})^{-1}.
\eqtn
Assume, in addition, that
\qtnl{180817z}
G_{ij}\ne G_{ik}\quad\text{for all}\ \,i,j,k\in I,\ j\ne k.
\eqtn
In this case if the indices $i,j,k\in I$ are distinct then each of the 
sections $S_{ijk}=G_i/G_{ij}G_{ik}$ and $S_{kij}=G_k/G_{ki}G_{kj}$ is trivial and hence formulas~\eqref{221216z} and~\eqref{050817a} hold. The same is obviously true if some of $i,j,k$ are equal. Therefore the triple
$$
\cT(\cG,\fF)=(\fG,\fS,\fF),
$$
where $\fS=\{G_i/G_{ij}\}_{i,j\in I}$ is a system of linked quotients. Denote by $\cX(\cG,\fF)$ the quasiregular coherent configuration associated with the triple $\cT(\cG,\fF)$.\medskip

Let us study the schurity question for a coherent configuration $\cX=\cX(\cG,\fF)$ for suitable families $\cG$ and $\fF$. Let~$a,i,j$ be distinct elements of~$I$. Then condition~\eqref{180817z} implies that
$$
G_{ia}\ne G_{ij}\qaq G_{ja}\ne G_{ji}
$$
and hence
$$
G_{ia}\times G_{ij}=G=G_{ja}\times G_{ji}.
$$ 
It follows that the isomorphism $f_{ija}:G_i/G_{ij}\to G_j/G_{ji}$ induces a unique isomorphism $\sigma_{ija}:G_{ia}\to G_{ja}$ such that
\qtnl{190817u}
(xG_{ij})^{f_{ij}}=x^{\sigma_{ija}}G_{ji},\qquad x\in G_{ia}.
\eqtn
Assume that the coherent configuration $\cX$ is schurian. For each $i\in I$, choose a point $\alpha_i\in\Omega_i$. Denote by $k$ the index other than $a,i,j$, i.e., $\{a,i,j,k\}=\{1,2,3,4\}$. Then the sets 
$$
\Omega_{ia}=\alpha_i G_{ia},\quad \Omega_{ja}=\alpha_jG_{ja},\quad\Omega_{ka}=\alpha_kG_{ka}
$$ 
are the neighborhoods of a point $\alpha_a$ in suitable basis relations of~$\cX$ (see Lemma~\ref{120817f}). It follows that these sets are orbits of the group $\aut(\cX)_\alpha$ for any $\alpha\in\Omega_a$. Moreover, for any two distinct indices $u,v\in\{i,j,k\}$, we have
$$
|r\,\cap\, (\Omega_{ua}\times \Omega_{va})|=p\quad\text{for all}\ \,r\in S_{uv},
$$
i.e., the relation $r_{uv}=r\,\cap\,(\Omega_{ua}\times \Omega_{va})$ is a matching. Note that $r_{uv}s_{vw}\in S_{uw}$ 
for all $r\in S_{uv}$ and $s\in S_{vw}$, where 
$w$ is the element of $\{i,j,k\}$ other than $u$ and~$v$. This immediately implies that
\qtnl{170617a}
\sigma_{ija}\sigma_{jka}\sigma_{kia}=\id_{G_i}.
\eqtn
This proves the following statement.

\prpstnl{170617b}
In the above notation, the coherent configuration $\cX$ is schurian only if formula~\eqref{170617a} holds true for all $i,j,k,a$ such that $\{i,j,k,a\}=\{1,2,3,4\}$.
\eprpstn

There are many ways how to choose the families $\cG$ and $\fF$ to satisfy the condition~\eqref{170617a}. To see a concrete example, let
$$
H_1=1\times C_p,\quad H_2=C_p\times 1,\quad H_3=\diag(C_p\times C_p).
$$
Denote by $\cG$ the family of the groups $G_{ij}$ that are $(i,j)$-entries
of the following $4\times 4$ array, $i,j\in I$:
$$
\begin{pmatrix}
1   & H_1 & H_2 & H_3\\
H_1 & 1   & H_3 & H_2\\
H_2 & H_3 & 1   & H_1\\
H_3 & H_2 & H_1 & 1\\
\end{pmatrix}.
$$
In particular, $G_{ij}=G_{ji}$ for all $i$ and~$j$. Next, each of the matrices 
$$
\sigma_{23}=\mattwo{0}{1}{1}{0},\quad
\sigma_{34}=\mattwo{1}{0}{1}{1},\quad
\sigma_{42}=\mattwo{1}{-1}{0}{1},
$$
with elements in $\mF_p$, belongs to $\GL(2,p)$ and hence induces an automorphism of the group~$G$ so that
$$
(H_1)^{\sigma_{23}}=H_2,\quad (H_2)^{\sigma_{34}}=H_3,\quad (H_3)^{\sigma_{42}}=H_1.
$$
Now, denote by $\fF$ the family the isomorphisms $f_{ij}:G/G_{ij}\to G/G_{ji}$ that are $(i,j)$-entries of the following $4\times 4$ array, $i,j\in I$:
$$
\begin{pmatrix}
\id_G   & \id_{G/H_1} & \id_{G/H_2} & \id_{G/H_3}\\
\id_{G/H_1}	& \id_G       & f_{23}      & f_{42}^{-1}\\
\id_{G/H_2}	& f_{23}^{-1} & \id_G       & f_{34}\\
\id_{G/H_3}	& f_{42}      & f_{34}^{-1} &  \id_G\\
\end{pmatrix},
$$
where the isomorphisms $f_{23}$, $f_{34}$, and $f_{42}$ are obtained from formula~\eqref{190817u} for $a=1$ and $\sigma_{ij1}=\sigma_{ij}|_{H_{i-1}}$ with $(i,j)=(2,3)$, $(3,4)$, and $(4,2)$, respectively. From the choice of the matrices $\sigma_{ij}$, it follows that 
$$
\sigma_{23}\sigma_{34}\sigma_{42}=\mattwo{1}{0}{1}{-1},
$$
which shows that condition~\eqref{170617a} is 
not satisfied for $(a,i,j,k)=(1,2,3,4)$. 
Thus, Proposition~\ref{170617b} implies the following statement.

\prpstnl{170817t}
For each prime~$p$, the coherent configuration $\cX(\cG,\fF)$ is not schurian. In particular, there are infinitely many non-schurian quasiregular coherent configurations with four fibers.
\eprpstn

We complete the section by a natural question, namely, whether the assumption in Proposition~\ref{170617b} gives a sufficient condition for the schurity of the coherent configuration $\cX(\cG,\fF)$? It would also be interesting to find an analog of that proposition for the separability problem.

\end{document}